\documentclass[12pt,a4paper]{article}
\usepackage{amsmath,amssymb,bm,ascmac,bbm,mathtools}
\usepackage[dvipdfmx, usenames]{color}
\usepackage[dvipdfmx]{graphicx}
\usepackage{xcolor}
\usepackage{here}
\usepackage{authblk}
\usepackage[hang,small,bf]{caption}
\usepackage[subrefformat=parens]{subcaption}
\usepackage{dcolumn}

\newcommand{\tr}{\text{tr}}

\newcommand{\mcal}{\mathcal}
\newcommand{\mbb}{\mathbb}
\newcommand{\mfrak}{\mathfrak}

\newcommand{\rank}{\text{rank}}

\newcommand{\mc}{\mathcal}

\newcommand{\wf}[1]{\widehat{\mfrak{#1}}}
\newcommand{\fp}{\text{FPdim}}
\newcommand{\fib}{\text{Fib}}
\newcommand{\ising}{\text{Ising}}
\newcommand{\vect}{\text{Vect}}
\newcommand{\vecG}{\text{Vec}}
\newcommand{\svec}{\text{sVec}}
\newcommand{\Hom}{\underline{\text{Hom}}}

\newcommand{\mods}{\text{mod }}

\begin{document}
\title{Classification of connected \'etale algebras in pre-modular fusion categories up to rank three}
\author{Ken KIKUCHI}
\affil{Department of Physics, National Taiwan University, Taipei 10617, Taiwan}
\date{}
\maketitle

\begin{abstract}
We classify connected étale algebras $A$'s in pre-modular fusion categories $\mathcal B$ with $\text{rank}(\mathcal B)\le3$ including degenerate and non-(pseudo-)unitary ones. We comment on Lagrangian algebras and physical applications to ground state degeneracy and proof of spontaneous $\mathcal B$-symmetry breaking.
\end{abstract}


\makeatletter
\renewcommand{\theequation}
{\arabic{section}.\arabic{equation}}
\@addtoreset{equation}{section}
\makeatother

\section{Introduction}
Throughout the paper, $\mcal C$ denotes a fusion category (FC) over the field $\mbb C$ of complex numbers (see \cite{EGNO15} for definitions). Its simple objects are denoted $c_i$'s with $i=1,\dots,\rank(\mcal C)$. We denote its fusion ring by $K(\mcal C)$. A pivotal structure $a:id_{\mcal C}\cong(-)^{**}$ defines left and right quantum (or categorical) traces $\tr^{L/R}(a):1\to1$. A pivotal structure is called spherical if $\forall c\in\mcal C$, $\tr^L(a_c)=\tr^R(a_{c^*})$. In an FC $\mcal C$ with a spherical structure $a$, a (quantum) dimension of an object $c\in\mcal C$ is defined as $d_c:=\tr^L(a_c)$. Its squared sum defines categorical (or global) dimension
\begin{equation}
    D^2(\mcal C):=\sum_{i=1}^{\rank(\mcal C)}d_{c_i}^2.\label{categoricaldim}
\end{equation}
Note that for a given categorical dimension, there are two $D(\mcal C)$'s, one positive and one negative. If an FC is equipped with a braiding $c$, in order to avoid confusion, we especially write the braided fusion category (BFC) $\mcal B$ and simple objects $b_i$'s. A spherical BFC is called pre-modular fusion category (pre-MFC). This is our ambient category. If the braiding in a pre-MFC is non-degenerate, the category is called modular or modular fusion category (MFC).

With a braiding, we can study commutative algebras. An algebra in an FC $\mcal C$ is a triple $(A,\mu,u)$ of an object $A\in\mcal C$, multiplication morphism $\mu:A\otimes A\to A$, and unit morphism $u:1\to A$ obeying associativity and unit axioms. A category $\mcal C_A$ of right $A$-modules consists of pairs $(m,p)$ where $m\in\mcal C$ and $p:m\otimes A\to m$ subject to coherence conditions. (A category $_A\mcal C$ of left $A$-modules is defined analogously.) An algebra is called separable if $\mcal C_A$ is semisimple. An algebra $A\in\mcal B$ in a BFC is called commutative if
\begin{equation}
    \mu\cdot c_{A,A}=\mu.\label{commutativealg}
\end{equation}
A commutative separable algebra is called étale. Any étale algebra $A$ decomposes to a direct sum of connected ones \cite{DMNO10} where $A\in\mcal C$ is called connected if $\dim_{\mbb C}\mcal C(1,A)=1$. An example of connected étale algebra is the unit $A\cong1\in\mcal B$ giving $\mcal B_A\simeq\mcal B$. A BFC $\mcal B$ without nontrivial connected étale algebra is called completely anisotropic. The goal of this paper is to classify connected étale algebras\footnote{They are also called quantum subgroup, consensable algebra, or normal algebra depending on contexts.} in pre-MFCs. Our main results are summarized as the\newline

\textbf{Theorem.} \textit{Connected étale algebras in pre-modular fusion categories $\mcal B$'s up to rank three are given as follows:}
\begin{table}[H]
\begin{center}
\begin{tabular}{c|c|c|c}
    Rank&$\mcal B$&Results&Completely anisotropic?\\\hline\hline
    1&$\vect_{\mbb C}$&Table \ref{rank1results}&Yes\\\hline
    2&$\vecG_{\mbb Z/2\mbb Z}^\alpha$&Tables \ref{rank2Z2anomfreeresults},\ref{rank2Z2anomresults}&No (two with (\ref{Z2etale}))/ Yes (the other 14)\\
    &$\fib$&Table \ref{rank2fibresults}&Yes\\\hline
    3&$\vecG_{\mbb Z/3\mbb Z}^1$&Table \ref{rank3Z3results}&No ($h_X=h_Y=0$)/ Yes ($h_X=h_Y=\pm\frac13$)\\
    &$\ising$&Table \ref{rank3isingresults}&Yes\\
    &$\text{Rep}(S_3)$&Tables \ref{rank3nonsymRepS3results},\ref{rank3symRepS3results}&No (the other six)/ Yes (two symmetric with $d_Y=-1$)\\
    &$psu(2)_5$&Table \ref{rank3psu25results}&Yes
\end{tabular}.
\end{center}
\caption{Connected étale algebras in pre-MFC $\mcal B$ with $\rank(\mcal B)\le3$}\label{results}
\end{table}

What is central in our classification method \cite{KK23GSD} is the Frobenius-Perron dimensions. They are defined as follows. An FC $\mcal C$ has Grothendieck ring $\text{Gr}(\mcal C)$. It admits a non-negative integer matrix representation (NIM-rep). The matrix $(N_i)_{jk}:={N_{ij}}^k$ consists of $\mbb N$-coefficients
\begin{equation}
    c_i\otimes c_j\cong\bigoplus_{k=1}^{\rank(\mcal C)}{N_{ij}}^kc_k.\label{NIM}
\end{equation}
Since they have non-negative entries, we can apply the Perron-Frobenius theorem to get the largest positive eigenvalue $\fp_{\mcal C}(c_i)$ called Frobenius-Perron dimension of $c_i$ (or $N_i$). Its squared sum defines Frobenius-Perron dimension of an FC
\begin{equation}
    \fp(\mcal C):=\sum_{i=1}^{\rank(\mcal C)}(\fp_{\mcal C}(c_i))^2.\label{FPdimC}
\end{equation}
An FC $\mcal C$ with $D^2(\mcal C)=\fp(\mcal C)$ is called pseudo-unitary, and unitary if $\forall c\in\mcal C$, $d_c=\fp_{\mcal C}(c)$. The Frobenius-Perron dimension defines a Lagrangian algebra. A connected étale algebra $A$ in a BFC $\mcal B$ is called Lagrangian if $(\fp_{\mcal B}(A))^2=\fp(\mcal B)$. With the Frobenius-Perron dimension at hand, our classification method \cite{KK23GSD} consists of three steps:
\begin{enumerate}
    \item Find a maximal rank $r_\text{max}$ on $\rank(\mcal B_A)$,
    \item List up candidate fusion categories for $\mcal B_A$,
    \item Check which of them satisfy axioms.
\end{enumerate}
A maximal rank on $\mcal B_A$ is given by $r_\text{max}=\lfloor\fp(\mcal B)\rfloor$ \cite{KK23GSD}. This is because $\forall c\in\mcal C$, $\fp_{\mcal C}(c)\ge1$ \cite{ENO02,EGNO15}, and for a BFC $\mcal B$ and a connected étale algebra $A\in\mcal B$, Frobenius-Perron dimensions obey \cite{KO01,ENO02,DMNO10}
\begin{equation}
    \fp_{\mcal B}(A)=\frac{\fp(\mcal B)}{\fp(\mcal B_A)}.\label{FPdimA}
\end{equation}
Since $\fp_{\mcal B}(A)\ge1$, we obtain
\begin{equation}
    1\le\fp(\mcal B_A)\le\fp(\mcal B),\label{fpdimeq}
\end{equation}
and
\begin{equation}
    \rank(\mcal B_A)\le r_\text{max}=\lfloor\fp(\mcal B)\rfloor.\label{rmax}
\end{equation}
In our setup, $\mcal B_A$ is known to be a fusion category.\footnote{Here is the reason. Let $\mcal C$ be a multitensor category and $A\in\mcal C$ be a connected separable algebra. Then, $\mcal C_A$ is a semisimple tensor category with simple unit $A\in\mcal C_A$. If $\mcal C$ is finite, $\mcal C_A$ is finite. Thus, for a finite multitensor category $\mcal C$ and connected separable algebra $A\in\mcal C$, $\mcal C_A$ is a fusion category. (Recall a fusion category is defined as finite semisimple tensor category.)} Hence, we list up candidate fusion categories $\mcal C$'s for $\mcal B_A$ in the second step. The candidates should satisfy three conditions: i) $\rank(\mcal C)\le r_\text{max}=\lfloor\fp(\mcal B)\rfloor$, ii) $(1\le)\fp(\mcal C)\le\fp(\mcal B)$, and iii) $\fp_{\mcal B}(A)=\fp(\mcal B)/\fp(\mcal C)$. The first two conditions are coming from (\ref{rmax},\ref{fpdimeq}). The third condition is coming from (\ref{FPdimA}). Since $A$ is a direct sum of simple objects in $\mcal B$, its Frobenius-Perron dimension is given by a linear sum of those of simple objects with $\mbb N$-coefficients
\begin{equation}
    \fp_{\mcal B}(A)=\sum_{i=1}^{\rank(\mcal B)}n_i\fp_{\mcal B}(b_i)\quad(n_i\in\mbb N).\label{fpdimAlinearsum}
\end{equation}
It turns out that these conditions are strong enough to rule out many fusion categories. In the third step, we check which candidate fusion categories satisfy axioms of connected étale algebras. In particular, we have to check (\ref{commutativealg}). For that purpose, we check a necessary condition
\begin{equation}
    \mu\cdot c_{A,A}\cdot c_{A,A}=\mu\label{necessarycommutative}
\end{equation}
for an algebra to be commutative. With conformal dimension $h_i$ of $b_i$, the double braiding is given by
\begin{equation}
    c_{A,A}\cdot c_{A,A}\cong\sum_{i,j=1}^n(\iota_i\otimes\iota_j)\cdot c_{a_j,a_i}\cdot c_{a_i,a_j}\cdot(p_i\otimes p_j),\label{cAAcAA}
\end{equation}
where $p_i:A\to a_i$ and $\iota_i:a_i\to A$ are product projections and coproduct injections of a direct sum $A\cong a_1\oplus a_2\oplus\cdots\oplus a_n$, respectively. The double braiding $c_{a_j,a_i}\cdot c_{a_i,a_j}$ of simple objects $a_i,a_j$ can be computed from the formula\footnote{The quantum trace of double braidings defines $S$-matrix
\begin{equation}
    \widetilde S_{i,j}:=\tr(c_{b_j,b_i}\cdot c_{b_i,b_j})=\sum_{k=1}^{\rank(\mcal B)}{N_{ij}}^k\frac{e^{2\pi ih_k}}{e^{2\pi i(h_i+h_j)}}d_k.\label{Smatrix}
\end{equation}
A normalized $S$-matrix is defined by
\[ S_{i,j}:=\frac{\widetilde S_{i,j}}{D(\mc B)}. \]
Another modular matrix $T$ is defined by
\begin{equation}
    T_{i,j}=e^{2\pi ih_i}\delta_{i,j}.\label{Tmatrix}
\end{equation}
When $\mcal B$ is modular, the matrices define (additive) central charge $c(\mcal B)$ mod 8 by
\begin{equation}
    (ST)^3=e^{2\pi ic(\mcal B)/8}S^2.\label{centralcharge}
\end{equation}}
\begin{equation}
    c_{a_j,a_i}\cdot c_{a_i,a_j}=\sum_{k=1}^{\rank(\mcal B)}{N_{ij}}^k\frac{e^{2\pi ih_k}}{e^{2\pi i(h_i+h_j)}}id_k.\label{doublebraiding}
\end{equation}
If one finds $A$ does not satisfy (\ref{necessarycommutative}), one can drop it from the list of candidates. If a candidate does pass the necessary condition, we need to check the axiom $\mu\cdot c_{A,A}=\mu$ seriously. Just like a double braiding, the half-braiding is given by
\begin{equation}
    c_{A,A}\cong\sum_{i,j=1}^n(\iota_j\otimes\iota_i)\cdot c_{a_i,a_j}\cdot(p_i\otimes p_j).\label{cAA}
\end{equation}
To find $c_{a_i,a_j}$, in principle, we have to solve hexagon equations, but `diagonal' ones $c_{b,b}$'s for self-dual simple objects $b\cong b^*$ can be computed from the formula \cite{K05}
\begin{equation}
    c_{b,b}\cong e^{-2\pi ih_b}\nu_2(b)id_{b\otimes b},\label{cbb}
\end{equation}
where $\nu_2(b)$ is the Frobenius-Schur indicator \cite{FRS92}.\footnote{When $\mcal B$ is modular, it can be computed with the formula
\begin{equation}
    \nu_2(b_k):=\frac1{D^2(\mcal B)}\sum_{i,j=1}^{\rank(\mcal B)}{N_{ij}}^kd_id_j\left(\frac{e^{2\pi ih_i}}{e^{2\pi ih_j}}\right)^2.\label{FSindicator}
\end{equation}}

Another fact we use to check an existence of connected étale algebra $A\in\mcal B$ is that $\mcal B_A$ is a left $\mcal B$-module category \cite{O01}. A left $\mcal C$-module category (or module category over $\mcal C$) is a quadruple $(\mcal M,\triangleright,m,l)$ of a category $\mcal M$, an action (or module product) bifunctor $\triangleright:\mcal C\times\mcal M\to\mcal M$, and natural isomorphisms $m_{-,-,-}:(-\otimes-)\triangleright-\cong-\triangleright(-\triangleright-)$ and $l:1\triangleright\mcal M\simeq\mcal M$ called module associativity constraint and unit constraint, respectively. They are subject to associativity and unit axioms as usual. An essential fact about module categories is that since $c\triangleright m$ is an object of $\mcal M$, it can be decomposed to a direct sum of simple objects with $\mbb N$-coefficients. Thus, actions of (left) $\mcal C$-module categories form NIM-reps. If one finds no $r$-dimensional NIM-rep, one can rule out candidates $A\in\mcal B$ giving $\mcal B_A$'s with rank $r$. On the other hand, if one finds $r$-dimensional NIM-reps, we can use it to constrain forms of $A$'s. In order to get additional constraints, we use the internal Hom \cite{O01}. For $m_1,m_2\in\mcal M$ over $\mcal C$, the internal Hom $\Hom(m_1,m_2)\in\mcal C$ is defined, if exists, by a natural isomorphism
\begin{equation}
    \forall c\in\mcal C,\quad\mcal M(c\triangleright m_1,m_2)\cong\mcal C(c,\Hom(m_1,m_2)).\label{inthom}
\end{equation}
It has a natural isomorphism
\begin{equation}
    \forall c\in\mcal C,\forall m_1,m_2\in\mcal M,\quad\Hom(m_1,c\triangleright m_2)\cong c\otimes\Hom(m_1,m_2).\label{inthomprop}
\end{equation}
We will use this fact to translate actions $\triangleright$ of $\mcal C$ on $\mcal M$ to monoidal products $\otimes$ in $\mcal C$. For every $m\in\mcal M$, $\Hom(m,m)$ has a canonical structure of an algebra in $\mcal C$. Furthermore, it is known \cite{O01,EGNO15} that the functor
\begin{equation}
    F:=\Hom(m,-):\mcal M\to\mcal C_{\Hom(m,m)}\quad(m\in\mcal M)\label{Ffunc}
\end{equation}
sending $m'\in\mcal M$ to $\Hom(m,m')\in\mcal C$ gives $\mcal M\simeq\mcal C_{\Hom(m,m)}$.

\section{Classification}
\subsection{Rank one}
A rank one pre-MFC has only one simple object $1$. There are two pre-MFCs depending on the sign of categorical dimension $D(\mcal B)=\pm1$. Since $1\le\fp(\mcal B_A)\le\fp(\mcal B)=1$, the only possibility is
\[ \fp(\mcal B_A)=1. \]
Accordingly, the algebra has
\[ \fp_{\mcal B}(A)=\frac{\fp(\mcal B)}{\fp(\mcal B_A)}=1. \]
Since the only simple object is $1\in\mcal B$, we get the unique algebra object
\begin{equation}
    A\cong1.\label{rank1alg}
\end{equation}
It is trivially connected étale. The algebra is also Lagrangian.

\begin{table}[H]
\begin{center}
\begin{tabular}{c|c|c|c}
    Connected étale algebra $A$&$\mcal B_A$&$\rank(\mcal B_A)$&Lagrangian?\\\hline
    1&$\mcal B$&1&Yes
\end{tabular}.
\end{center}
\caption{Connected étale algebra in rank one pre-MFC $\mcal B\simeq\vect_{\mbb C}$}\label{rank1results}
\end{table}
\hspace{-17pt}

\subsection{Rank two}
A rank two pre-MFC has two simple objects $\{1,X\}$. There are two fusion rings, $\mbb Z/2\mbb Z$ and $K(\fib)$. We start from the $\mbb Z/2\mbb Z$ case.

\subsubsection{$\mcal B\simeq\vecG_{\mbb Z/2\mbb Z}^\alpha$}
The pre-MFCs $\mcal B$'s have two simple objects $\{1,X_h\}$. They obey the $\mbb Z/2\mbb Z$ fusion ring
\[ \forall b\in\mcal B,\quad1\otimes b\cong b\cong b\otimes1,\quad X\otimes X\cong1, \]
and have
\[ \fp_{\mcal B}(1)=1=\fp_{\mcal B}(X). \]
Accordingly, the pre-MFCs have
\[ \fp(\mcal B)=2. \]
The quantum dimension $d_X$ is a solution of $d_X^2=1$, and there are two possible values
\[ d_X=\pm1. \]
They both have the same categorical dimension
\[ D^2(\mcal B)=2. \]
The $\mbb Z/2\mbb Z$ object $X$ can have conformal dimensions
\[ h=0,\frac14,\frac12,\frac34\quad(\mods1). \]
Those with $h=0,1/2$ mod 1 are anomaly-free $(\alpha=1)$, and the others $h=1/4,3/4$ mod 1 are anomalous $(\alpha=-1)$. Therefore, for the $\mbb Z/2\mbb Z$ fusion ring, we have
\[ 4(\text{conformal dimension }h)\times2(\text{quantum dimension }d_X)\times2(\text{categorical dimension }D=\pm\sqrt2)=16 \]
pre-MFCs. For a pre-MFC to be unitary, we need $d_X=1$. Furthermore, for such pre-MFCs to be modular, we need $h=\frac14,\frac34$.\footnote{A pre-MFC with $X_h$ has $S$-matrix (in the basis $\{1,X\}$)
\[ \widetilde S=\begin{pmatrix}1&d_X\\d_X&e^{-4\pi ih}\end{pmatrix}. \]} Therefore, there are
\[ 2(\text{conformal dimension }h)\times1(d_X=1)\times2(\text{categorical dimension }D)=4 \]
unitary MFCs, consistent with \cite{RSW07}. However, we do not assume non-degeneracy nor (pseudo-)unitarity. Rather, we search for connected étale algebras in all 16 pre-MFCs simultaneously.

First, since $\fp(\mcal B)=2$, an upper bound on ranks of $\mcal B_A$'s is
\begin{equation}
    r_\text{max}=2.\label{rank2Z2rmax}
\end{equation}
Fusion categories up to rank two are completely classified \cite{O02}, and there are only three fusion categories, $\vect_{\mbb C}$, $\vecG_{\mbb Z/2\mbb Z}^\alpha$, and $\fib$. The Fibonacci fusion ring has $\fp=\frac{5+\sqrt5}2>\fp(\mcal B)$, and is ruled out. The other two satisfy three conditions, and remain as candidates.

As the third step, we check which fusion categories satisfy axioms. We start from the candidate $\vecG_{\mbb Z/2\mbb Z}^\alpha$. The fusion category has $\fp=2$, and (\ref{FPdimA}) demands $\fp_{\mcal B}(A)=1$. For such an object to be connected, the only candidate is
\begin{equation}
    A\cong1.\label{rank2Z2rank2alg}
\end{equation}
This is the trivial connected étale algebra giving $\mcal B_A\simeq\mcal B$. Next, we study the rank one candidate $\vect_{\mbb C}$. The fusion category has $\fp=1$, and (\ref{FPdimA}) demands $\fp_{\mcal B}(A)=2$. For such an object to be connected, the only candidate is
\begin{equation}
    A\cong1\oplus X.\label{rank2Z2rank1alg}
\end{equation}
The candidate is connected, however, fails to be commutative depending on conformal dimensions. Let us study when the candidate is commutative. Since $X^*\cong X$ is invertible, we have \cite{K05}
\begin{equation}
    \alpha=F^{XXX}_X=\frac{\nu_2(X)}{d_X}.\label{FSindX}
\end{equation}
The formula (\ref{cbb}) gives
\[c_{X,X}\cong e^{-2\pi ih_X}\alpha d_Xid_1\cong\begin{cases}id_1&(d_X,h)=(1,0),(-1,\frac12),\\-id_1&(d_X,h)=(1,\frac12),(-1,0),\\+i\cdot id_1&(d_X,h)=(1,\frac14),(-1,\frac34),\\-i\cdot id_1&(d_X,h)=(1,\frac34),(-1,\frac14).\end{cases} \]
Therefore, the candidate (\ref{rank2Z2rank1alg}) is commutative iff $(d_X,h)=(1,0),(-1,\frac12)$ mod 1 for $h$. In other words, our pre-MFC $\mcal B\simeq\vecG_{\mbb Z/2\mbb Z}^\alpha$ can admit nontrivial connected étale algebra if it is anomaly-free $\alpha=1$, and anomalous $\alpha=-1$ ones are completely anisotropic. The result has a simple (partial) physical explanation. When the symmetry $\mbb Z/2\mbb Z$ is anomaly-free, we can gauge it. Thus, it is expected that there are more algebra objects. On the other hand, if the symmetry is anomalous, it cannot be gauged and there cannot exist a $\mbb Z/2\mbb Z$ algebra object. (This is the reason why there are fewer possible values of ground state degeneracies in anomalous theories compared to anomaly-free theories.)

Note that not all commutative algebras are separable. For $A\cong1\oplus X$ to be separable, it should have quantum dimension one. In $\mc B_A$, we have \cite{KO01}
\[ d_{\mc B_A}(c)=\frac{d_{\mc B}(c)}{d_{\mc B}(A)}, \]
where $d_{\mc C}(c)$ is the quantum dimension of $c$ as an object of $\mc C$. When $d_X=-1$, this cannot be one. Therefore, we conclude
\begin{equation}
    A\cong1\oplus X\quad(d_X,h_X)=(1,0).\quad(\mods1\text{ for }h)\label{Z2etale}
\end{equation}

To summarize, we arrive
\begin{table}[H]
\begin{center}
\begin{tabular}{c|c|c|c}
    Connected étale algebra $A$&$\mcal B_A$&$\rank(\mcal B_A)$&Lagrangian?\\\hline
    1&$\mcal B$&2&No\\
    $1\oplus X$ for (\ref{Z2etale})&$\vect_{\mbb C}$&1&No
\end{tabular},
\end{center}
\caption{Connected étale algebra in rank two pre-MFC $\mcal B\simeq\vecG_{\mbb Z/2\mbb Z}^1$}\label{rank2Z2anomfreeresults}
\end{table}
\hspace{-17pt}and
\begin{table}[H]
\begin{center}
\begin{tabular}{c|c|c|c}
    Connected étale algebra $A$&$\mcal B_A$&$\rank(\mcal B_A)$&Lagrangian?\\\hline
    1&$\mcal B$&2&No
\end{tabular}.
\end{center}
\caption{Connected étale algebra in rank two pre-MFC $\mcal B\simeq\vecG_{\mbb Z/2\mbb Z}^{-1}$}\label{rank2Z2anomresults}
\end{table}

\subsubsection{$\mcal B\simeq\fib$}
The pre-MFCs $\mcal B$'s have two simple objects $\{1,X_h\}$ obeying the Fibonacci fusion ring
\[ \forall b\in\mcal B,\quad1\otimes b\cong b\cong b\otimes1,\quad X\otimes X\cong1\oplus X. \]
They thus have
\[ \fp_{\mcal B}(1)=1,\quad\fp_{\mcal B}(X)=\zeta:=\frac{1+\sqrt5}2, \]
and
\[ \fp(\mcal B)=\frac{5+\sqrt5}2\approx3.6. \]
The quantum dimension $d_X$ is a solution of $d_X^2-d_X-1=0$, and it can take two values
\[ d_X=\zeta,-\zeta^{-1}. \]
Accordingly, they have categorical dimensions
\[ D^2(\mcal B)=\frac{5+\sqrt5}2,\quad\frac{5-\sqrt5}2, \]
respectively, and modular matrices
\[ \widetilde S=\begin{pmatrix}1&\zeta\\\zeta&-1\end{pmatrix},\quad T=\begin{pmatrix}1&\\&e^{2\pi ih}\end{pmatrix}\quad(d_X=\zeta), \]
or
\[ \widetilde S=\begin{pmatrix}1&-\zeta^{-1}\\-\zeta^{-1}&-1\end{pmatrix},\quad T=\begin{pmatrix}1&\\&e^{2\pi ih}\end{pmatrix}\quad(d_X=-\zeta^{-1}). \]
The former has
\[ h=\frac25,\frac35\quad(\mods1), \]
and the latter has
\[ h=\frac15,\frac45\quad(\mods1). \]
Therefore, we have
\[ 2(\text{quantum dimension }d_X)\times2(\text{conformal dimension }h)\times2(\text{categorical dimension }D(\mcal B))=8 \]
pre-MFCs. Those with $d_X=\zeta$ are unitary, and there are four unitary MFCs, consistent with \cite{RSW07}. We study connected étale algebras in all eight pre-MFCs simultaneously.

First, since $\fp(\mcal B)\approx3.6$, an upper bound on ranks of $\mcal B_A$'s is
\begin{equation}
    r_\text{max}=3.\label{rank2fibrmax}
\end{equation}
Fusion categories up to rank three have been classified in \cite{O02,O13} assuming an existence of pivotal structure for rank three. There are
\[ 1(\text{rank}=1)+2(\text{rank}=2)+5(\text{rank}=3)=8 \]
fusion rings. All but one are free of multiplicity, and they are summarized in AnyonWiki \cite{anyonwiki} based on \cite{LPR20,VS22}. The one with multiplicity has fusion ring $K(1,0,2,0)$ in the notation of \cite{O13}, however, it has Frobenius-Perron dimension\footnote{Simple objects $\{1,X,Y\}$ obey fusion rules
\begin{table}[H]
\begin{center}
\begin{tabular}{c|c|c|c}
    $\otimes$&1&$X$&$Y$\\\hline
    $1$&1&$X$&$Y$\\\hline
    $X$&&$1\oplus2X\oplus Y$&$X$\\\hline
    $Y$&&&$1$
\end{tabular},
\end{center}
\end{table}
\hspace{-17pt}
and have
\[ \fp_{\mcal C}(1)=1=\fp_{\mcal C}(Y),\quad\fp_{\mcal C}(X)=1+\sqrt3. \]} $\fp(\mcal C)=6+2\sqrt3\approx9.4$, and we can rule it out at the outset. As in \cite{KK23GSD}, one can easily check there are no one- and three-dimensional NIM-reps. Hence, we are left with rank two $\mcal B_A$'s. They have fusion rings\footnote{The notation in \cite{anyonwiki} is the following. A rank $r$ fusion ring with $n$ non-self-dual simple objects is denoted by $\text{FR}^{r,n}_i$. The subscript $i=1,\dots$ labels different fusion rings with the same $r$ and $n$.}
\[ \text{FR}^{2,0}_1,\text{FR}^{2,0}_2. \]
In order to search for all rank two $\mcal B_A$'s, we solve two-dimensional NIM-reps. We find a unique (up to basis transformations) solution
\[ n_1=1_2,\quad n_X=\begin{pmatrix}0&1\\1&1\end{pmatrix}. \]
Let $m_1,m_2\in\mcal M$ be basis. The NIM-rep gives actions
\begin{align*}
    1\triangleright m_j&\cong m_j,\\
    X\triangleright m_1&\cong m_2,\\
    X\triangleright m_2&\cong m_1\oplus m_2,
\end{align*}
or monoidal products
\begin{align*}
    F(m_j)\cong F(1\triangleright m_j)\equiv\Hom(m,1\triangleright m_j)&\cong\Hom(m,m_j),\\
    F(m_2)\cong F(X\triangleright m_1)\equiv\Hom(m,X\triangleright m_1)&\cong X\otimes\Hom(m,m_1),\\
    F(m_1)\oplus F(m_2)\cong F(X\triangleright m_2)\equiv\Hom(m,X\triangleright m_2)&\cong X\otimes\Hom(m,m_2).
\end{align*}
In order to find forms of algebras, we set an ansatz
\[ A\cong F(m)\cong a1\oplus bX, \]
with $a,b\in\mbb N$. If $m\cong m_1$, the second monoidal product yields
\[ F(m_2)\cong b1\oplus(a+b)X. \]
The third monoidal product is automatically satisfied. For the candidate to be connected, we need $a=1$, or $A\cong F(m_1)\cong1\oplus bX$. It has $\fp_{\mcal B}(A)=1+b\zeta$. For this to solve (\ref{FPdimA}), we need $b=0$ and fusion ring $\text{FR}^{2,0}_2$. Namely, we only get the trivial connected étale algebra $A\cong1$ giving $\mcal B_A\simeq\mcal B$. The other case $m\cong m_2$ gives the same result.

In conclusion, we arrive
\begin{table}[H]
\begin{center}
\begin{tabular}{c|c|c|c}
    Connected étale algebra $A$&$\mcal B_A$&$\rank(\mcal B_A)$&Lagrangian?\\\hline
    1&$\mcal B$&2&No
\end{tabular}.
\end{center}
\caption{Connected étale algebra in rank two pre-MFC $\mcal B\simeq\fib$}\label{rank2fibresults}
\end{table}
\hspace{-17pt}Note that since rank three $\mcal B_A$'s are ruled out, we do not have to assume $\mcal B_A$ admits a pivotal structure or its fusion ring be multiplicity-free. The result implies all $\fib$'s are completely anisotropic, giving another proof of the known fact \cite{BD11}.

\subsection{Rank three}
A rank three pre-MFC has three simple objects $\{1,X,Y\}$. As we reviewed above, with the assumption of pivotal structure, there are five fusion rings, $\text{FR}^{3,2}_1,\text{FR}^{3,0}_1,\text{FR}^{3,0}_2,\text{FR}^{3,0}_3,K(1,0,2,0)$. The last with multiplicity is known \cite{HH07} not to admit braiding. Thus, we classify connected étale algebras in the first four\footnote{The fourth with fusion ring $\text{FR}^{3,0}_3$ does not seem to have a common name. Thus, we call it $psu(2)_5$ borrowing one famous realization.} pre-MFCs.

\subsubsection{$\mcal B\simeq\vecG_{\mbb Z/3\mbb Z}^1$}\label{Z3}
The pre-MFCs $\mcal B$'s have three simple objects $\{1,X_{h_X},Y_{h_Y}\}$ obeying $\mbb Z/3\mbb Z$ fusion ring
\begin{table}[H]
\begin{center}
\begin{tabular}{c|c|c|c}
    $\otimes$&1&$X$&$Y$\\\hline
    $1$&1&$X$&$Y$\\\hline
    $X$&&$Y$&$1$\\\hline
    $Y$&&&$X$
\end{tabular}.
\end{center}
\end{table}
\hspace{-17pt}Hence, they have
\[ \fp_{\mcal B}(1)=\fp_{\mcal B}(X)=\fp_{\mcal B}(Y)=1, \]
and
\[ \fp(\mcal B)=3. \]
The quantum dimensions $d_X,d_Y$ are solutions of $d_X^2=d_Y,\ d_Y^2=d_X,\ d_Xd_Y=1$. The only solution is
\[ d_X=d_Y=1. \]
Thus, the pre-MFCs are unitary. Consequently, the pre-MFCs have a unique categorical dimension
\[ D^2(\mcal B)=3. \]
The nontrivial simple objects are dual to each other, $X^*\cong Y$, and they have the same conformal dimensions
\[ h_X=h_Y=0,\frac13,\frac23\quad(\mods1). \]
This implies the pre-MFCs are always anomaly-free $\alpha=1$. (We have already specified this fact in the name of this subsubsection.) Therefore, there are
\[ 3(\text{conformal dimension }h_X=h_Y)\times2(\text{categorical dimension }D(\mcal B)=\pm\sqrt3)=6 \]
pre-MFCs with $\mbb Z/3\mbb Z$ fusion ring. Among them, those with nontrivial conformal dimensions are modular. There are four of them, consistent with \cite{RSW07}. We study all six pre-MFCs including degenerate ones simultaneously.

In order to classify connected étale algebras in our pre-MFCs, we first find an upper bound on $\rank(\mcal B_A)$. It is given by
\begin{equation}
    r_\text{max}=3.\label{rank3Z3rmax}
\end{equation}
As we reviewed above, there are eight fusion categories which admit pivotal structures, but the condition $\fp(\mcal B_A)\le\fp(\mcal B)=3$ narrows down to three with fusion rings
\[ \text{FR}^{1,0}_1,\text{FR}^{2,0}_1,\text{FR}^{3,2}_1. \]
The rank two fusion ring cannot solve (\ref{FPdimA}), and is ruled out. This leaves only two candidates, $\text{FR}^{1,0}_1,\text{FR}^{3,2}_1$. Let us study the two candidates in detail.

We start from the rank three candidate. It has $\fp(\mcal C)=3$, and demands $\fp_{\mcal B}(A)=1$. This is nothing but the trivial connected étale algebra $A\cong1$ giving $\mcal B_A\simeq\mcal B$. Next, we study the rank one candidate $\vect_{\mbb C}$. It has $\fp=1$, and demands $\fp_{\mcal B}(A)=3$. Since any étale algebras are self-dual \cite{DMNO10}, we are left with two possibilities, $1\oplus1\oplus1$ or $1\oplus X\oplus Y$. Connectedness of algebra rules out the first. Thus, the only candidate is $A\cong1\oplus X\oplus Y$. This is connected and separable, but it fails to be commutative depending on conformal dimensions. Let us first show those with nontrivial conformal dimensions are non-commutative. We check the necessary condition (\ref{necessarycommutative}). The relevant double-braidings are given by
\[ c_{X,X}\cdot c_{X,X}\cong e^{-2\pi ih_X}id_Y,\quad c_{Y,X}\cdot c_{X,Y}\cong e^{-4\pi ih_X}id_1,\quad c_{Y,Y}\cdot c_{Y,Y}\cong e^{-2\pi ih_Y}id_X. \]
(We used $h_X=h_Y$.) The phases are nontrivial for nontrivial conformal dimensions $h_X=h_Y=\frac13,\frac23$ mod 1, and they fail to satisfy the necessary condition. Hence, they are non-commutative. On the other hand, when the conformal dimensions are trivial $h_X=h_Y=0$ mod 1, the phases become trivial, and they pass the necessary condition. We thus need to compute half-braidings $c_{b,b'}$'s. We find\footnote{Employing the known \cite{RSW07} associators $F^{b_i,b_j,b_k}_{b_l}=1$, one can easily solve hexagon equations to find $R$-matrices
\[ R^{X,X}_Y=R^{X,Y}_1=R^{Y,X}_1=R^{Y,Y}_X=+1, \]
giving the half-braidings.}
\[ c_{X,X}\cong id_Y,\quad c_{X,Y}\cong id_1,\quad c_{Y,Y}\cong id_X. \]
Therefore, the candidate $A\cong1\oplus X\oplus Y$ is commutative for trivial conformal dimensions.

To summarize, we found
\begin{table}[H]
\begin{center}
\begin{tabular}{c|c|c|c}
    Connected étale algebra $A$&$\mcal B_A$&$\rank(\mcal B_A)$&Lagrangian?\\\hline
    1&$\mcal B$&3&No\\
    $1\oplus X\oplus Y$ for $h_X=h_Y=0$ (mod 1)&$\vect_{\mbb C}$&1&No
\end{tabular}.
\end{center}
\caption{Connected étale algebra in rank three pre-MFC $\mcal B\simeq\vecG_{\mbb Z/3\mbb Z}^1$}\label{rank3Z3results}
\end{table}
\hspace{-17pt}In other words, the pre-MFCs $\vecG_{\mbb Z/3\mbb Z}^1$'s with nontrivial conformal dimensions $h_X=h_Y=\frac13,\frac23$ mod 1 are completely anisotropic.

\subsubsection{$\mcal B\simeq\text{Ising}$}
The pre-MFCs $\mcal B$'s have three simple objects $\{1,X_{h_X},Y_{h_Y}\}$ obeying the Ising fusion ring
\begin{table}[H]
\begin{center}
\begin{tabular}{c|c|c|c}
    $\otimes$&1&$X$&$Y$\\\hline
    $1$&1&$X$&$Y$\\\hline
    $X$&&$1$&$Y$\\\hline
    $Y$&&&$1\oplus X$
\end{tabular}.
\end{center}
\end{table}
\hspace{-17pt}Hence, they have
\[ \fp_{\mcal B}(1)=1=\fp_{\mcal B}(X),\quad\fp_{\mcal B}(Y)=\sqrt2, \]
and
\[ \fp(\mcal B)=4. \]
It is known \cite{DGNO09} that Ising fusion categories admit spherical structures and (non-degenerate) braidings. (We will make the most of the non-degeneracy later.) The quantum dimensions $d_X,d_Y$ are solutions of $d_X^2=1,\ d_Xd_Y=d_Y,\ d_Y^2=1+d_X$, and there are two solutions
\[ d_X=1,\quad d_Y=\pm\sqrt2. \]
(The sign of $d_Y$ fixes the spherical structure.) The pre-MFCs with positive $d_Y$ are unitary. The pre-MFCs have a unique categorical dimension
\[ D^2(\mcal B)=4. \]
It is known that conformal dimensions take values
\[ h_X=\frac12,\quad h_Y=\pm\frac1{16},\ \pm\frac3{16},\ \pm\frac5{16},\ \pm\frac7{16}\quad(\mods1). \]
(The conformal dimension $h_X=\frac12$ implies the rank two braided fusion subcategory generated by $\{1,X\}$ is isomorphic to $\svec$. Hence, an obvious $\mbb Z/2\mbb Z$ algebra $A\cong1\oplus X$ fails to be commutative because $c_{X,X}\cong-id_1$.) Therefore, there are
\[ 2(\text{quantum dimension }d_Y)\times8(\text{conformal dimension }h_Y)\times2(\text{categorical dimension }D(\mcal B))=32 \]
pre-MFCs. Those with $d_Y=+\sqrt2$ are unitary, and there are 16 of them, consistent with \cite{RSW07}. We search for connected étale algebras in all 32 pre-MFCs simultaneously. Since the method is the same, we will be brief.

The Frobenius-Perron dimension $\fp(\mcal B)=4$ gives an upper bound
\begin{equation}
    r_\text{max}=4.\label{rank4Isingrmax}
\end{equation}
Since there are no complete list of fusion categories up to rank four, here we practically have to assume fusion rings of $\mcal B_A$'s be multiplicity-free. Then, there are seven\footnote{We could take the rank three fusion ring with multiplicity into account, but it is ruled out because $\fp>\fp(\mcal B)=4$.} fusion rings, but only five
\[ \text{FR}^{1,0}_1,\text{FR}^{2,0}_1,\text{FR}^{3,0}_1,\text{FR}^{4,0}_1,\text{FR}^{4,2}_1 \]
satisfy the three conditions. Let us check which of them satisfy axioms.

We start from those with $\fp=4$. The relation (\ref{FPdimA}) requires $\fp_{\mcal B}(A)=1$. For such an algebra to be connected, the only possibility is the trivial one $A\cong1$ giving $\mcal B_A\simeq\mcal B$. (This rules out rank four $\mcal B_A$'s.) Next, let us see the rank one scenario. One finds there is no one-dimensional NIM-rep, and the fusion ring $\text{FR}^{1,0}_1$ is ruled out. Finally, let us study the last candidate $\text{FR}^{2,0}_1$. The scenario demands $\fp_{\mcal B}(A)=2$. For such an algebra to be connected, the only possibility is $A\cong1\oplus X$. Indeed, this is an algebra,\footnote{An easiest way to see this fact is to realize that $X$ generates anomaly-free $\mbb Z/2\mbb Z$ symmetry in physical systems.} however, it fails to be commutative because $X$ is a fermion $c_{X,X}\cong-id_1$. Thus, we do not get new connected étale algebra from the rank two scenario.

To summarize, we found
\begin{table}[H]
\begin{center}
\begin{tabular}{c|c|c|c}
    Connected étale algebra $A$&$\mcal B_A$&$\rank(\mcal B_A)$&Lagrangian?\\\hline
    1&$\mcal B$&3&No
\end{tabular}
\end{center}
\caption{Connected étale algebra in rank three pre-MFC $\mcal B\simeq\text{Ising}$}\label{rank3isingresults}
\end{table}
\hspace{-17pt}regardless of quantum dimension $d_Y$, conformal dimension $h_Y$, or categorical dimension $D(\mcal B)$. The result motivates the\newline

\textbf{Theorem.} \textit{All }$\mcal B\simeq\ising$\textit{'s are completely anisotropic.}\newline

\textit{Proof.} So far, we pretended we do not know $\mcal B\simeq\ising$ be modular. However, as we recalled above, they are known to be modular \cite{DGNO09}. Thus, we can employ further constraints. When the ambient category $\mcal B$ is modular, it is known \cite{KO01,DMNO10} that the category of dyslectic modules $\mcal B_A^0$ for connected étale algebra $A$ is modular obeying
\begin{equation}
\begin{split}
    \fp(\mcal B_A^0)&=\frac{\fp(\mcal B)}{(\fp_{\mcal B}(A))^2},\\
    e^{2\pi ic(\mcal B)/8}&=e^{2\pi ic(\mcal B_A^0)/8}.
\end{split}\label{constsMFCB}
\end{equation}
Here, the full subcategory $\mcal B_A^0$ of $\mcal B_A$ consists of dyslectic modules $(m,p)\in\mcal B_A$ such that $p\cdot c_{A,m}\cdot c_{m,A}=p$. Since $\mcal B_A^0\subset\mcal B_A$, we also have
\[ \rank(\mcal B_A^0)\le\rank(\mcal B)\le\lfloor\fp(\mcal B)\rfloor. \]
In our example, we learn the maximal rank of $\mcal B_A^0$ is four. Fortunately, MFCs up to rank four are completely classified \cite{GK94,RSW07}.\footnote{See also \cite{BNRW15} for classification of MFCs at rank five.} According to the lists, the only MFCs which can match the (additive) central charges are Ising MFCs. Since the MFCs have $\fp=4$, it demands $\mcal B_A^0\simeq\mcal B_A\simeq\mcal B$, and $\fp_{\mcal B}(A)=1$. Therefore, the only connected étale algebra is the trivial one $A\cong1$, showing $\mcal B\simeq\ising$ be completely anisotropic. $\square$\newline

\textbf{Remark.} We have already included the theorem in the Table \ref{rank3isingresults}.\newline

\textbf{Remark.} Some Ising MFCs are realized by Wess-Zumino-Witten (WZW) models:
\begin{table}[H]
\begin{center}
\begin{tabular}{c|c}
    Unitary $\ising$ w/ $h_Y$ (mod 1)&Realization\\\hline
    $3/16$&$\wf{su}(2)_2$\\\hline
    $5/16$&$\wf{sp}(4)_1$\\\hline
    $7/16$&$\wf{so}(7)_1$\\\hline
    $15/16$&$(\wf{e_8})_2$
\end{tabular}.
\end{center}
\end{table}
\hspace{-17pt}Connected étale algebras in MFCs $\mcal C(\mfrak g,k)$ describing $\wf g_k$ WZW models are especially called quantum subgroups \cite{O00}. They have been studied extensively, and many results are known. Such connected étale algebras were classified in \cite{KO01} (for $\wf{su}(2)_k$), in \cite{EP09} (for $\wf{su}(3)_k$), in \cite{CEM23} (for $\wf{su}(4)_k$), and many more in \cite{G23} (see also \cite{EM21}). Thus, the classification results for some Ising MFCs are known. Our results are consistent with them.\footnote{For $\wf{su}(2)_2$, $\mcal B_A$ corresponds to the $A_3$ Dynkin diagram. Thus, the category $\mcal B_A$ has three simple objects, implying $A\cong1$. For $\wf{sp}(4)_1$ and $\wf{so}(7)_1$, they do not admit exotic quantum subgroups. That is, all nontrivial algebras are their centers $\mbb Z/2\mbb Z$. However, simple objects generating them are fermions, and the $\mbb Z/2\mbb Z$ algebras fail to be commutative. Thus, they are also completely anisotropic. For $(\wf{e_8})_2$, its only nontrivial algebra is the $\mbb Z/2\mbb Z$ \cite{G23}, but it is fermionic and fails to be commutative. Thus, the MFC is also completely anisotropic, and our results are consistent.

Of course, the critical Ising model is described by an Ising MFC with $h_Y=1/16$. The MFC is known (see, e.g., \cite{FRS03}) to be completely anisotropic, and our results are consistent. (The MFC has only one nontrivial connected separable algebra $A\cong1\oplus X$, but since $X$ is a fermion $c_{X,X}\cong-id_1$, it fails to be commutative.) See also \cite{KL02} for classification of connected étale algebras in minimal models.} However, we could not find a literature classifying connected étale algebras in the other Ising MFCs (in particular, non-unitary ones). Our results for those MFCs would be new. This example demonstrates a unified nature of our method; since what is central in our method is the Frobenius-Perron dimensions, it allows us to classify connected étale algebras in pre-MFCs (or more generally BFCs) with the same fusion ring -- including degenerate and non-(pseudo-)unitary ones -- simultaneously.

\subsubsection{$\mcal B\simeq\text{Rep}(S_3)$}
The pre-MFCs $\mcal B$'s have three simple objects $\{1,X_{h_X},Y_{h_Y}\}$ obeying the fusion ring
\begin{table}[H]
\begin{center}
\begin{tabular}{c|c|c|c}
    $\otimes$&1&$X$&$Y$\\\hline
    $1$&1&$X$&$Y$\\\hline
    $X$&&$1$&$Y$\\\hline
    $Y$&&&$1\oplus X\oplus Y$
\end{tabular}.
\end{center}
\end{table}
\hspace{-17pt}Hence, they have
\[ \fp_{\mcal B}(1)=1=\fp_{\mcal B}(X),\quad\fp_{\mcal B}(Y)=2, \]
and
\[ \fp(\mcal B)=6. \]
The quantum dimensions $d_X,d_Y$ are (nonzero) solutions to $d_X^2=1,\ d_Xd_Y=d_Y,\ d_Y^2=1+d_X+d_Y$, and there are two of them
\[ d_X=1,\quad d_Y=-1,2. \]
Accordingly, there are two categorical dimensions
\[ D^2(\mcal B)=3,6. \]
As we reviewed in Appendix \ref{repS3confdim}, there are two classes of pre-MFCs with fusion ring $K(\text{Rep}(S_3))$, one non-symmetric, and another symmetric. The non-symmetric pre-MFCs have conformal dimensions
\[ h_X=0,\quad h_Y=\frac13,\frac23\quad(\mods1) \]
with $d_Y=2$, and symmetric pre-MFCs have
\[ h_X=h_Y=0\quad(\mods1). \]
Therefore, there are
\[ 2(\text{sign of categorical dimension }D(\mcal B))\times2(h_Y=\pm\frac13)=4 \]
non-symmetric and
\[ 2(\text{sign of categorical dimension }D(\mcal B))\times2(d_Y=-1,2)=4 \]
symmetric pre-MFCs, eight in total. Since none of them is modular, they should not be listed in MFCs, consistent with \cite{RSW07}. We classify connected étale algebras in all eight pre-MFCs simultaneously.

First, since $\fp(\mcal B)=6$, the rank of categories of right $A$-modules is bounded from above by
\begin{equation}
    r_\text{max}=6.\label{rank3RepD3rmax}
\end{equation}
Fusion categories up to rank six have been classified in \cite{LPR20,VS22} assuming multiplicity-free fusion rings. In order to apply our method, for a moment, we assume $\mcal B_A$ have multiplicity-free fusion ring. Later, employing an improved method, we relax the assumption and complete the classification of connected étale algebras.

According to \cite{anyonwiki}, we find seven fusion rings
\[ \text{FR}^{1,0}_1,\text{FR}^{2,0}_1,\text{FR}^{3,2}_1,\text{FR}^{3,0}_2,\text{FR}^{4,2}_2,\text{FR}^{6,2}_1,\text{FR}^{6,4}_1 \]
can describe $\mcal B_A$'s. The same analysis was done in \cite{KK23GSD}, and we found only four
\[ \text{FR}^{1,0}_1,\text{FR}^{2,0}_1,\text{FR}^{3,2}_1,\text{FR}^{3,0}_2 \]
can satisfy axioms. They correspond to candidates
\[ A\cong1\oplus X\oplus2Y,\quad1\oplus Y,\quad1\oplus X,\quad1, \]
respectively. Which of them are connected étale? The answer depends on conformal dimensions, or non/symmetry of our ambient category. We start from the non-symmetric ones.

\paragraph{Non-symmetric $\text{Rep}(S_3)$.} In this case, the simple object $Y$ has nontrivial braiding. The first two candidates fail to satisfy the necessary condition (\ref{necessarycommutative}),\footnote{They have double braidings
\[ c_{A,A}\cdot c_{A,A}\cong(2+4e^{-4\pi ih_Y})\iota\cdot id_1\cdot p+(2+e^{-4\pi ih_Y})\iota\cdot id_X\cdot p+(8+4e^{-2\pi ih_Y})\iota\cdot id_Y\cdot p, \]
and
\[ c_{A,A}\cdot c_{A,A}\cong(1+e^{-4\pi ih_Y})\iota\cdot id_1\cdot p+e^{-4\pi ih_Y}\iota\cdot id_X\cdot p+(2+e^{-2\pi ih_Y})\iota\cdot id_Y\cdot p, \]
respectively.} and we can safely discard them. We do know the last two are connected étale algebras.

\paragraph{Symmetric $\text{Rep}(S_3)$.} In this case, the simple object $Y$ also has trivial conformal dimension, and the candidates containing $Y$ can be commutative. We have to check whether the first two satisfy axioms of commutative algebra. Naively, one has to solve associativity axioms, however, luckily, we do not need any computation.\footnote{In order to check the associativity axioms directly, one may use the method developed in \cite{FRS02} and associators spelled out in, say, \cite{BBW21}. We checked the solution in \cite{CZW18} is consistent.} The second candidate $A\cong1\oplus Y$ is known \cite{CZW18} to be connected étale when $d_Y=2$.\footnote{The condition is needed because if $d_Y=-1$, the candidate has
\[ d_{\mc B}(A)=d_{\mc B}(1)+d_{\mc B}(Y)=0, \]
a contradiction.} (They call it condensable or normal.) Thus, we are left with the first candidate $A\cong1\oplus X\oplus2Y$. Note that we are in a special situation; our ambient category is symmetric, and if $A$ defines an algebra, it gives $\mcal B_A\simeq\vect_{\mbb C}$. Since the condensation gives a tensor functor $F:\mcal B\to\mcal B_A$, our question can be rephrased as follows. Is a symmetric fusion category $\mcal B\simeq\text{Rep}(S_3)$ Tannakian? (A symmetric fusion category admitting a braided tensor functor $F:\mcal B\to\vect_{\mbb C}$ is called Tannakian. Since our $\mcal B$ is symmetric and the only simple object of $\vect_{\mbb C}$ is the unit, a tensor functor is automatically braided.) Once we have rewritten our problem in this form, we can immediately answer it; $\mcal B\simeq\text{Rep}(S_3)$ is Tannakian iff it is positive \cite{D02,EGNO15}. Recall a symmetric fusion category $\mcal B$ is called positive if $\forall b\in\mcal B$, $d_b\in\mbb N^\times$.
Namely, $A\cong1\oplus X\oplus2Y$ is connected étale\footnote{The half-braidings are given by
\[ c_{X,X}\cong id_1,\quad c_{X,Y}\cong id_Y,\quad c_{Y,Y}\cong id_1+id_X+id_Y. \]
One of the easiest ways to get these is as follows. First, since we know $A\cong1\oplus Y$ is commutative, we learn
\[ c_{Y,Y}\cong id_1+id_X+id_Y. \]
Then, one can easily solve hexagon equations to get the results above. From the half-braidings, one immediately finds
\[ c_{A,A}\cong6\Big(\iota\cdot id_1\cdot p+\iota\cdot id_X\cdot p+2\iota\cdot id_Y\cdot p\Big)\cong id_{A\otimes A}. \]} if $d_Y=2$, but not otherwise. Finally, for $A\cong1\oplus X$ to be separable, we also need $d_Y=2$. The reason is as follows. In view of anyon condensation, $Y$ `splits' into two simple objects. Thus, simple objects of $\mc B_A$ are given by
\[ 1\oplus X,\quad Y,\quad Y. \]
They have quantum dimensions \cite{KO01} $d_{\mc B_A}(c)=\frac{d_{\mc B}(c)}{d_{\mc B}(A)}$. Since the only rank three fusion category with Frobenius-Perron dimension three is $\vecG_{\mbb Z/3\mbb Z}^1$, for $A$ to be separable, the three simple objects should form $\vecG_{\mbb Z/3\mbb Z}^1$. As we saw in section \ref{Z3}, simple objects of the fusion category have quantum dimensions one. Thus, for $A$ to give a fusion category, the simple objects need to have quantum dimensions one. This condition is satisfied only when $d_Y=2$.
\newline

So far, we assumed $K(\mc B_A)$ be multiplicity-free. However, employing the improved method \cite{KK23MFCrank9}, we can relax the assumption.

We work with an ansatz
\[ A\cong1\oplus n_XX\oplus n_YY \]
with $n_j\in\mbb N$. It has
\[ \fp_{\mc B}(A)=1+n_X+2n_Y. \]
For this to obey
\[ 1\le\fp_{\mc B}(A)\le\fp(\mc B) \]
resulting from (\ref{FPdimA}), the natural numbers can take only 12 values
\begin{align*}
    (n_X,n_Y)=&(0,0),(1,0),(0,1),(2,0),(1,1),(0,2),\\
    &(3,0),(2,1),(1,2),(4,0),(3,1),(5,0).
\end{align*}
The first is nothing but the trivial connected étale algebra $A\cong1$ giving $\mc B_A^0\simeq\mc B_A\simeq\mc B$. Among the other 11 candidates, five are ruled out by studying Frobenius-Perron dimensions. The two candidates $(n_X,n_Y)=(1,1),(3,0)$ have $\fp_{\mc B}=4$, and demands $\fp(\mc B_A)=\frac32$. However, there is no such fusion category, and the candidates are ruled out. Similarly, three candidates $(n_X,n_Y)=(0,2),(2,1),(4,0)$ have $\fp_{\mc B}=5$, and demand $\fp(\mc B_A)=\frac65$, but there is no such fusion category. Thus, the three candidates are also ruled out. We are left with six nontrivial candidates:
\[ (n_X,n_Y)=(1,0),(0,1),(2,0),(1,2),(3,1),(5,0). \]
Regardless of conformal dimensions, three of them, $(n_X,n_Y)=(2,0),(3,1),(5,0)$ or $A\cong1\oplus2X,1\oplus3X\oplus Y,1\oplus5X$, are ruled out by detailed study.

\paragraph{$A\cong1\oplus2X$.} It has $\fp_{\mc B}(A)=3$, and demands $\fp(\mc B_A)=2$. With the free module functor $F_A:=-\otimes A:\mc B\to\mc B_A$, we find candidate simple objects of $\mc B_A$:
\[ 1\oplus2X,\quad2\oplus X,\quad3Y. \]
They have Frobenius-Perron dimensions $1,1,2$, respectively. Their contributions to $\fp(\mc B_A)$ exceed two, and the candidate is ruled out.

\paragraph{$A\cong1\oplus3X\oplus Y$.} It has $\fp_{\mc B}(A)=6$, and demands $\fp(\mc B_A)=1$. The candidate simple objects of $\mc B_A$ are given by
\[ 1\oplus3X\oplus Y,\quad3\oplus X\oplus Y,\quad\dots\ . \]
There exists additional simple objects, but their presence does not affect the following discussion. Their contributions to $\fp(\mc B_A)$ exceed one, and the candidate is ruled out.

\paragraph{$A\cong1\oplus5X$.} It has $\fp_{\mc B}(A)=6$, and demands $\fp(\mc B_A)=1$. We obtain candidate simple objects (assuming the smallest Frobenius-Perron dimensions)
\[ 1\oplus5X,\quad5\oplus X,\quad3Y. \]
They have $\fp_{\mc B_A}=1,1,1$, and their contributions to $\fp(\mc B_A)$ exceed one. Thus, the candidate is also ruled out.

Now, we are only left with three nontrivial candidates
\[ (n_X,n_Y)=(1,0),(0,1),(1,2). \]
We have seen that these give nontrivial connected étale algebras depending on quantum and conformal dimensions. From this analysis, we succeeded to relax the assumption of multiplicity on $K(\mc B_A)$.

To summarize, we found
\begin{table}[H]
\begin{center}
\begin{tabular}{c|c|c|c}
    Connected étale algebra $A$&$\mcal B_A$&$\rank(\mcal B_A)$&Lagrangian?\\\hline
    1&$\mcal B$&3&No\\
    $1\oplus X$&$\vecG_{\mbb Z/3\mbb Z}^1$&3&No
\end{tabular}
\end{center}
\caption{Connected étale algebra in rank three non-symmetric pre-MFC $\mcal B\simeq\text{Rep}(S_3)$}\label{rank3nonsymRepS3results}
\end{table}
\hspace{-17pt}for non-symmetric, and
\begin{table}[H]
\begin{center}
\begin{tabular}{c|c|c|c}
    Connected étale algebra $A$&$\mcal B_A$&$\rank(\mcal B_A)$&Lagrangian?\\\hline
    1&$\mcal B$&3&No\\
    $1\oplus X$ for $d_Y=2$&$\vecG_{\mbb Z/3\mbb Z}^1$&3&No\\
    $1\oplus Y$ for $d_Y=2$&$\vecG_{\mbb Z/2\mbb Z}^\alpha$&2&No\\
    $1\oplus X\oplus2Y$ for $d_Y=2$&$\vect_{\mbb C}$&1&No
\end{tabular}
\end{center}
\caption{Connected étale algebra in rank three symmetric pre-MFC $\mcal B\simeq\text{Rep}(S_3)$}\label{rank3symRepS3results}
\end{table}
\hspace{-17pt}for symmetric $\mc B\simeq\text{Rep}(S_3)$. Thus, two symmetric ones with $d_Y=-1$ are completely anisotropic, while the other six $\mcal B\simeq\text{Rep}(S_3)$'s fail to be completely anisotropic.

\subsubsection{$\mcal B\simeq psu(2)_5$}
The pre-MFCs $\mcal B$'s have three simple objects $\{1,X_{h_X},Y_{h_Y}\}$ obeying the fusion ring $\text{FR}^{3,0}_3$
\begin{table}[H]
\begin{center}
\begin{tabular}{c|c|c|c}
    $\otimes$&1&$X$&$Y$\\\hline
    $1$&1&$X$&$Y$\\\hline
    $X$&&$1\oplus Y$&$X\oplus Y$\\\hline
    $Y$&&&$1\oplus X\oplus Y$
\end{tabular}.
\end{center}
\end{table}
\hspace{-17pt}Hence, they have
\[ \fp_{\mcal B}(1)=1,\quad\fp_{\mcal B}(X)=\frac{\sin\frac{2\pi}7}{\sin\frac\pi7}\approx1.802,\quad\fp_{\mcal B}(Y)=\frac{\sin\frac{3\pi}7}{\sin\frac\pi7}\approx2.247, \]
and
\[ \fp(\mcal B)=\frac7{4\sin^2\frac\pi7}\approx9.3. \]
The quantum dimensions $d_X,d_Y$ are solutions of $d_X^2=1+d_Y,\ d_Y^2=1+d_X+d_Y,\ d_Xd_Y=d_X+d_Y$. There are three solutions
\[ (d_X,d_Y)=(\frac{\sin\frac\pi7}{\cos\frac\pi{14}},-\frac{\sin\frac{2\pi}7}{\cos\frac\pi{14}}),\ (-\frac{\sin\frac{3\pi}7}{\cos\frac{3\pi}{14}},\frac{\sin\frac\pi7}{\cos\frac{3\pi}{14}}),\ (\frac{\sin\frac{2\pi}7}{\sin\frac\pi7},\frac{\sin\frac{3\pi}7}{\sin\frac\pi7}), \]
and categorical dimensions
\[ D^2(\mcal B)=\frac7{4\cos^2\frac\pi{14}},\quad\frac7{4\cos^2\frac{3\pi}{14}},\quad\frac7{4\sin^2\frac\pi7}, \]
respectively. The pre-MFCs have conformal dimensions
\begin{equation}
    (h_X,h_Y)=\begin{cases}(\frac37,\frac17),(\frac47,\frac67)&(\text{1st quantum dimensions}),\\
    (\frac27,\frac37),(\frac57,\frac47)&(\text{2nd quantum dimensions}),\\
    (\frac17,\frac57),(\frac67,\frac27)&(\text{3rd quantum dimensions}),\end{cases}\quad(\mods1)\label{psu25confdim}
\end{equation}
respectively. Therefore, there are
\[ 3(\text{quantum dimension})\times2(\text{sign of categorical dimension }D(\mcal B))\times2(\text{conformal dimension})=12 \]
pre-MFCs. (Actually, they are all modular. See the lemma in Appendix \ref{psu25}.) Those with the third quantum dimensions are unitary, and there are four of them, consistent with \cite{RSW07}. We study connected étale algebras in all 12 MFCs.

As the first step, we find an upper bound
\begin{equation}
    r_\text{max}=9\label{rank3psu25rmax}
\end{equation}
on $\rank(\mcal B_A)$.\footnote{Actually, we can do a little more. We can show a\newline

\textbf{Lemma.} \textit{An upper bound on $\rank(\mcal B_A^0)$ is eight. Namely, $\rank(\mcal B_A^0)$ cannot be nine.}\newline

\textit{Proof.} Assume the opposite, and suppose $\rank(\mcal B_A^0)=9$. Then, it is either pointed or non-pointed. If it is pointed, we only have two possibilities, $\mcal B_A^0\simeq\vecG_{\mbb Z/9\mbb Z}^1,\vecG_{\mbb Z/3\mbb Z\times\mbb Z/3\mbb Z}^1$. The MFCs with such fusion rings cannot match the additive central charges. Thus, we are left with non-pointed possibilities. A rank nine non-pointed fusion category $\mcal C$ with $\fp(\mcal C)\le\fp(\mcal B)=\frac7{4\sin^2\frac\pi7}$ is strongly constrained. Since it is non-pointed, there exists a non-invertible object $c\in\mcal C$ with $\fp_{\mcal C}(c)>1$. It enters the Frobenius-Perron dimension by squared sum, $\fp(\mcal C)\equiv\sum_{i=1}^9(\fp_{\mcal C}(c_i))^2=(\fp_{\mcal C}(c))^2+\cdots$. Since any object has $\fp_{\mcal C}(c_i)\ge1$, $\sum_{c_i\neq c}(\fp_{\mcal C}(c_i))^2\ge8$, and we must have $1<\fp_{\mcal C}(c)<2$, or the squared sum reaches ten. Here, it is known \cite{ENO02} that such Frobenius-Perron dimension is expressed as $\fp_{\mcal C}=2\cos\frac\pi n$ for a natural number $n$ larger than three ($n=3$ is ruled out by our assumption $\fp_{\mcal C}>1$). Since the function is monotonically increasing, the smallest Frobenius-Perron dimension is achieved when $n=4$, or $\fp_{\mcal C}=\sqrt2$. The Frobenius-Perron dimensions of such fusion categories exceed the upper bound $\fp(\mcal B)=\frac7{4\sin^2\frac\pi7}\approx9.3$, leading to a contradiction. $\square$} According to \cite{anyonwiki}, we find $\mcal B_A$ can have only two fusion rings (assuming $K(\mcal B_A)$ be multiplicity-free)
\[ \text{FR}^{1,0}_1,\text{FR}^{3,0}_3. \]
The second scenario requires $\fp_{\mcal B}(A)=1$. For such a candidate to be connected, the only candidate is $A\cong1$. This is the trivial connected étale algebra giving $\mcal B_A\simeq\mcal B$. To check the first candidate with rank one, we search for one-dimensional NIM-reps. One finds no solution, ruling out the possibility of rank one $\mcal B_A$.

To summarize, we found
\begin{table}[H]
\begin{center}
\begin{tabular}{c|c|c|c}
    Connected étale algebra $A$&$\mcal B_A$&$\rank(\mcal B_A)$&Lagrangian?\\\hline
    1&$\mcal B$&3&No
\end{tabular}.
\end{center}
\caption{Connected étale algebra in rank three pre-MFC $\mcal B\simeq psu(2)_5$}\label{rank3psu25results}
\end{table}
\hspace{-17pt}The result motivates the\newline

\textbf{Theorem.} \textit{All $\mcal B\simeq psu(2)_5$'s are completely anisotropic.}\newline

\textit{Proof.} As we recalled above, all pre-MFCs are actually modular. Thus, the algebra should obey (\ref{constsMFCB}). Since $\fp(\mcal B_A^0)\ge1$, we obtain an inequality\footnote{We thank Victor Ostrik for teaching us an upper bound on $\fp_{\mcal B}(A)$ in the related work \cite{KK23GSD}.}
\begin{equation}
    1\le\left(\fp_{\mcal B}(A)\right)^2\le\fp(\mcal B).\label{FPdimA2bound}
\end{equation}
The most general candidate can be written as
\[ A\cong\bigoplus_{j=1}^{\rank(\mcal B)}n_jb_j \]
with $n_j\in\mbb N$. In order to get connected algebra, we set $n_j=1$ for $b_j\cong1$ at the outset. The candidate in our example thus has
\[ \fp_{\mcal B}(A)=1+\frac1{\sin\frac\pi7}\left(n_X\sin\frac{2\pi}7+n_Y\sin\frac{3\pi}7\right). \]
For this to obey (\ref{FPdimA2bound}), there are only two possibilities
\[ (n_X,n_Y)=(0,0),(1,0). \]
The first candidate is nothing but the trivial connected étale algebra $A\cong1$ giving $\mcal B_A^0\simeq\mcal B_A\simeq\mcal B$. The second candidate means $A\cong1\oplus X$. It has double braiding
\[ c_{A,A}\cdot c_{A,A}\cong(1+e^{-4\pi ih_X})\iota\cdot id_1\cdot p+2\iota\cdot id_X\cdot p+e^{2\pi i(h_Y-2h_X)}\iota\cdot id_Y\cdot p. \]
Since $X$ has nontrivial conformal dimensions, the candidate fails to satisfy the necessary condition (\ref{necessarycommutative}), and it is ruled out. $\square$
\newline

\textbf{Remark.} We have already included the theorem in the Table \ref{rank3psu25results}.

\section{Physical applications}
In this section, we discuss physical applications of the classification results.

Consider a two-dimensional $\mcal C$-symmetric gapped phase (an existence of braiding is not assumed at this level). It is known \cite{TW19,HLS21} that such phases stand in bijection with $\mcal C$-module categories
\begin{equation}
    \{\text{2d }\mcal C\text{-symmetric gapped phases}\}\cong\{\mcal C\text{-module categories }\mcal M\}.\label{1to1}
\end{equation}
In particular, ground state degeneracy (GSD) in the LHS is given by $\rank(\mcal M)$ of a module category $\mcal M$ in the RHS. Therefore, we can translate physical problems in the LHS to mathematical problems in the RHS. One of the most important physical problems is whether the $\mcal C$ symmetry is spontaneously broken or not. In this setup, a spontaneous symmetry breaking (SSB) in the LHS is defined as follows:\newline

\textbf{Definition.} \cite{KK23GSD} Let $\mcal C$ be a fusion category and $\mcal M$ be a (left) $\mcal C$-module category describing a $\mcal C$-symmetric gapped phase. A symmetry $c\in\mcal C$ is called \textit{spontaneously broken} if $\exists m\in\mcal M$ such that $c\triangleright m\not\cong m$. We also say $\mcal C$\textit{ is spontaneously broken} if there exists a spontaneously broken object $c\in\mcal C$. A categorical symmetry $\mcal C$ is called \textit{preserved} (i.e., not spontaneously broken) if all objects act trivially.\newline

With the definition, we can show a\newline

\textbf{Lemma.} \cite{KK23GSD} \textit{Let $\mcal C$ be a fusion category and $\mcal M$ be an indecomposable (left) $\mcal C$-module category. Then, $\rank(\mcal M)>1$ implies SSB of $\mcal C$ (i.e., $\mcal C$ is spontaneously broken).}\newline

Therefore, in two-dimensional gapped phases, we can mathematically prove SSB of a (categorical) symmetry $\mcal C$ by classifying $\mcal C$-module categories. In particular, if $\mcal M$ is finite, it is known \cite{O01,EGNO15} that $\exists A\in\mcal C$ such that $\mcal M\simeq\mcal C_A$. In classifying such module categories, it is natural to assume $A$ be connected étale. (Note that we also assume ambient categories be pre-modular at this point. There are infinitely many examples in this class, relevant deformations of rational conformal field theories.) Then, our classification results can be translated to physics side. We arrive the\newline

\textbf{Theorem.} \textit{Let $\mcal B$ be a pre-modular fusion category and $A\in\mcal B$ a connected étale algebra. Suppose two-dimensional $\mcal B$-symmetric gapped phases are described by indecomposable $\mcal B_A$'s. Then, the gapped phases have}
\[ \text{GSD}\in\begin{cases}\{1\}&(\mcal B\simeq\vect_{\mbb C}),\\\{1,2\}&(\mcal B\simeq\vecG_{\mbb Z/2\mbb Z}^1\ (\text{two with }(\ref{Z2etale}))),\\\{2\}&(\mcal B\simeq\vecG_{\mbb Z/2\mbb Z}^\alpha\text{ (the other 14)}),\\\{2\}&(\mcal B\simeq\fib),\\\{1,3\}&(\mcal B\simeq\vecG_{\mbb Z/3\mbb Z}^1\text{ with }h_X=h_Y=0),\\\{3\}&(\mcal B\simeq\vecG_{\mbb Z/3\mbb Z}^1\text{ with }h_X=h_Y=\pm\frac13),\\\{3\}&(\mcal B\simeq\ising),\\\{1,2,3\}&(\mcal B\simeq\text{Rep}(S_3)\ (\text{two with }(d_Y,h_Y)=(2,0)),\\\{3\}&(\mcal B\simeq\text{Rep}(S_3)\ (\text{the other six})),\\\{3\}&(\mcal B\simeq psu(2)_5).\end{cases} \]
\textit{In particular, categorical symmetries $\mcal B$'s are spontaneously broken for}
\[ \mcal B\simeq\begin{cases}\vecG_{\mbb Z/2\mbb Z}^\alpha\text{ with }(d_X,h_X)\neq(1,0),\\\fib,\\\vecG_{\mbb Z/3\mbb Z}^1\text{ with }h_X=h_Y=\pm\frac13,\\\ising,\\\text{Non-symmetric }\text{Rep}(S_3),\\\text{Symmetric }\text{Rep}(S_3)\text{ with }d_Y=-1,\\psu(2)_5.\end{cases} \]\newline

\textbf{Remark.} As noted in \cite{KK23GSD} Example 6, commutativity of an algebra seems too strong; a numerical computation suggests an existence of two-dimensional $\mcal B$-symmetric gapped phase described by $\mcal B_A$ with non-commutative connected separable algebra $A\in\mcal B$.


\appendix
\setcounter{section}{0}
\renewcommand{\thesection}{\Alph{section}}
\setcounter{equation}{0}
\renewcommand{\theequation}{\Alph{section}.\arabic{equation}}

\section{Conformal dimensions of $\mcal B\simeq\text{Rep}(S_3)$}\label{repS3confdim}
There are two classes of pre-MFCs with fusion ring $\text{FR}^{3,0}_2$ (or $K(0,1,0,1)$ in the notation of \cite{O05}); both are degenerate, but one is symmetric and the other is non-symmetric. In this appendix, we compute conformal dimensions of them.

\paragraph{Symmetric $\text{Rep}(S_3)$.} The symmetric $\text{Rep}(S_3)$ has $S$-matrix (in the basis $\{1,X,Y\}$)
\[ \widetilde S=\begin{pmatrix}1&1&d_Y\\1&1&d_Y\\d_Y&d_Y&d_Y^2\end{pmatrix}. \]
(We used $d_X=1$.) Since the $(i,j)$ component is defined by $\widetilde S_{i,j}:=\tr(c_{b_j,b_i}\cdot c_{b_i,b_j})=\sum_{k=1}^{\rank(\mcal B)}{N_{i,j}}^k\frac{e^{2\pi ih_k}}{e^{2\pi i(h_i+h_j)}}d_k$, we can read off the conformal dimensions. The $(2,3)$ component says $e^{2\pi ih_X}=1$, or
\begin{equation}
    h_X=0\quad(\mods1).\label{RepS3hX}
\end{equation}
Similarly, the $(3,3)$ component imposes $d_Y^2=2e^{-4\pi ih_Y}+e^{-2\pi ih_Y}d_Y$, or
\begin{equation}
    h_Y=0\quad(\mods1).\label{RepS3hY}
\end{equation}

\paragraph{Non-symmetric $\text{Rep}(S_3)$.} The non-symmetric $\text{Rep}(S_3)$ has $S$-matrix \cite{O05}
\[ \widetilde S=\begin{pmatrix}1&1&d_Y\\1&1&d_Y\\d_Y&d_Y&e^{-4\pi ih_Y}(2+e^{2\pi ih_Y}d_Y)\end{pmatrix} \]
with $(3,3)$ component equals to $(-2)$. (We already used $e^{2\pi ih_X}=1$, or $h_X=0$ mod 1.) Solving
\[ -2=e^{-4\pi ih_Y}(2+e^{2\pi ih_Y}d_Y) \]
for $d_Y=2,-1$, we find
\begin{equation}
    h_X=0,\quad h_Y=\begin{cases}\frac13,\frac23&(d_Y=2),\\\text{no solution}&(d_Y=-1),\end{cases}\quad(\mods1).\label{RepS3nonsymconfdim}
\end{equation}
The possibility $h_Y=\frac13$ is realized by even primaries in the $\wf{su}(2)_4$ WZW model, and another $h_Y=\frac23$ is realized by rank three pre-modular fusion subcategory of MFC describing $M(6,5)$ minimal model. Thus, both of them exist. For $d_Y=-1$, we naively get $h_Y=\frac1{2\pi}\arccos\frac14$, however, this contradicts Vafa's theorem \cite{V88}; $\exists N\in\mbb N^\times$ such that $\left(e^{2\pi ih_Y}\right)^N=1$. Thus, we conclude non-symmetric pre-MFC with fusion ring $K(\text{Rep}(S_3))$ must have $d_Y=2$ and $h_Y=\pm\frac13$ mod 1.

\section{Conformal dimensions of $\mcal B\simeq psu(2)_5$}\label{psu25}
By definition, the $S$-matrix in this example is given by (in the basis $\{1,X,Y\}$)
\[ \widetilde S=\begin{pmatrix}1&d_X&d_Y\\d_X&e^{-4\pi ih_X}+e^{2\pi i(h_Y-2h_X)}d_Y&e^{-2\pi ih_Y}d_X+e^{-2\pi ih_X}d_Y\\d_Y&e^{-2\pi ih_Y}d_X+e^{-2\pi ih_X}d_Y&e^{-4\pi ih_Y}+e^{2\pi i(h_X-2h_Y)}d_X+e^{-2\pi ih_Y}d_Y\end{pmatrix}. \]
The $S$-matrix obeys $\widetilde S_{i^*,j}=\widetilde S_{i,j}^*$. Since all simple objects are self-dual $b^*\cong b$, this means all components should be real. This imposes three necessary conditions on conformal dimensions $h_X,h_Y$:
\begin{equation}
\begin{split}
    e^{-4\pi ih_X}+e^{2\pi i(h_Y-2h_X)}d_Y&\stackrel!\in\mbb R,\\
    e^{-2\pi ih_Y}d_X+e^{-2\pi ih_X}d_Y&\stackrel!\in\mbb R,\\
    e^{-4\pi ih_Y}+e^{2\pi i(h_X-2h_Y)}d_X+e^{-2\pi ih_Y}d_Y&\stackrel!\in\mbb R.
\end{split}\label{psu25realSmatrixcond}
\end{equation}
Solving the conditions for $(d_X,d_Y)=(\frac{\sin\frac\pi7}{\cos\frac\pi{14}},-\frac{\sin\frac{2\pi}7}{\cos\frac\pi{14}}),\ (-\frac{\sin\frac{3\pi}7}{\cos\frac{3\pi}{14}},\frac{\sin\frac\pi7}{\cos\frac{3\pi}{14}}),\ (\frac{\sin\frac{2\pi}7}{\sin\frac\pi7},\frac{\sin\frac{3\pi}7}{\sin\frac\pi7})$, one finds

\begin{equation}
    (h_X,h_Y)=\begin{cases}(\frac37,\frac17),(\frac47,\frac67)&(\text{1st quantum dimensions}),\\
    (\frac27,\frac37),(\frac57,\frac47)&(\text{2nd quantum dimensions}),\\
    (\frac17,\frac57),(\frac67,\frac27)&(\text{3rd quantum dimensions}).\end{cases}\label{psu25confdimresult}
\end{equation}
Indeed, these conformal dimensions give known forms \cite{RSW07} of $S$-matrices, and the MFCs exist.

Naively, we also find trivial conformal dimensions $h_X=0=h_Y$ (mod 1) satisfy the necessary conditions. However, we can rule them out by a\newline

\textbf{Lemma.} \textit{There does not exist a pre-modular fusion category with fusion ring $K(psu(2)_5)$ and trivial conformal dimension. In other words, all pre-modular fusion categories with fusion ring $K(psu(2)_5)$ are modular.}\newline

\textit{Proof.} Assume the opposite, and suppose there exist pre-MFCs with trivial conformal dimensions obeying the fusion ring $K(psu(2)_5)$. When conformal dimensions are trivial, they turn out to be symmetric. For a symmetric fusion category $\mcal B$, there exists a finite group $G$ such that $\mcal B\simeq\text{Rep}(G)$ \cite{D02}. The only possibilities in our case is $G=\mbb Z/3\mbb Z,S_3$ \cite{O05}, and such symmetric fusion categories should have fusion rings $K(\text{Rep}(\mbb Z/3\mbb Z))$ or $K(\text{Rep}(S_3))=K(0,1,0,1)$ (in the notation of the paper). However, the fusion ring of $psu(2)_5$ is different from neither of them, leading to a contradiction. $\square$


\begin{thebibliography}{30}
\bibitem{EGNO15}
  P.~Etingof, S.~Gelaki, D.~Nikshych and V.~Ostrik, ``Tensor Categories,'' American Mathematical Society, 2015.
\bibitem{DMNO10}
  A.~Davydov, M.~Müger, D.~Nikshych and V.~Ostrik, ``The Witt group of non-degenerate braided fusion categories,'' Journal für die reine und angewandte Mathematik (Crelles Journal), 2013(677), 135-177. https://doi.org/10.1515/crelle.2012.014
  [arXiv:1009.2117 [math.QA]].
\bibitem{KK23GSD}
  K.~Kikuchi,
  ``Ground state degeneracy and module category,''
  [arXiv:2311.00746 [hep-th]].
\bibitem{ENO02}
  P.~Etingof, D.~Nikshych and V.~Ostrik,
  ``On fusion categories,''
  Annals of Mathematics 162, no. 2 (2005): 581–642. http://www.jstor.org/stable/20159926
  [arXiv:math/0203060 [math.QA]].
\bibitem{KO01}
  A.~Kirillov Jr. and V.~Ostrik, ``ON A q-ANALOG OF THE MCKAY CORRESPONDENCE AND THE ADE CLASSIFICATION OF slb2 CONFORMAL FIELD,'' Advances in Mathematics 171(2002), 183-227. https://doi.org/10.1006/aima.2002.2072 [arXiv:math/0101219 [math.QA]].
\bibitem{K05}
  A.~Kitaev,
  ``Anyons in an exactly solved model and beyond,''
  Annals Phys. \textbf{321}, no.1, 2-111 (2006)
  doi:10.1016/j.aop.2005.10.005
  [arXiv:cond-mat/0506438 [cond-mat.mes-hall]].
\bibitem{FRS92}
  K.~Fredenhagen, K.~H.~Rehren and B.~Schroer,
  ``Superselection sectors with braid group statistics and exchange algebras. 2. Geometric aspects and conformal covariance,''
  Rev. Math. Phys. \textbf{4}, no.spec01, 113-157 (1992)
  doi:10.1142/S0129055X92000170
\bibitem{O01}
  V.~Ostrik, ``MODULE CATEGORIES, WEAK HOPF ALGEBRAS AND MODULAR INVARIANTS,'' Transformation Groups \textbf8, 177–206 (2003). https://doi.org/10.1007/s00031-003-0515-6
  [arXiv:math/0111139 [math.QA]].
\bibitem{RSW07}
  E.~Rowell, R.~Stong, Z.~Wang,
  ``On classification of modular tensor categories,'' Communications in Mathematical Physics 292(2009), 343-389. https://doi.org/10.1007/s00220-009-0908-z
  [arXiv:0712.1377 [math.QA]].
\bibitem{O02}
  V.~Ostrik, ``Fusion categories of rank 2,'' Mathematical Research Letters 10 (2002): 177-183. https://dx.doi.org/10.4310/MRL.2003.v10.n2.a5
  [arXiv:math/0203255 [math.QA]].
\bibitem{O13}
  V.~Ostrik, ``Pivotal fusion categories of rank 3,''  Mosc. Math. J., 15(2015), 373–396. https://doi.org/10.17323/1609-4514-2015-15-2-373-396
  [arXiv:1309.4822 [math.QA]].
\bibitem{anyonwiki}
  ``AnyonWiki,'' https://anyonwiki.github.io/
\bibitem{LPR20}
  Z.~Liu, S.~Palcoux and Y.~Ren, ``Classification of Grothendieck rings of complex fusion categories of multiplicity one up to rank six,'' Lett Math Phys 112, 54 (2022). https://doi.org/10.1007/s11005-022-01542-1
  [arXiv:2010.10264 [math.CT]].
\bibitem{VS22}
  G.~Vercleyen and J.~Slingerland, ``On Low Rank Fusion Rings,''
  [arXiv:2205.15637 [math-ph]].
\bibitem{BD11}
  T.~Booker and A.~Davydov, ``Commutative Algebras in Fibonacci Categories,'' Journal of Algebra 355(2012), 176-204. https://doi.org/10.1016/j.jalgebra.2011.12.029
  [arXiv:1103.3537 [math.CT]].
\bibitem{HH07}
  T.~J.~Hagge and S.-M.~Hong, ``Some non-braided fusion categories of rank 3,'' Communications in Contemporary Mathematics 11(2009), 615-637. https://doi.org/10.1142/S0219199709003521
  [arXiv:0704.0208 [math.GT]].
\bibitem{DGNO09}
  V.~Drinfeld, S.~Gelaki, D.~Nikshych and V.~Ostrik, ``On braided fusion categories I,'' Selecta Mathematica 16(2010), 1-119. https://doi.org/10.1007/s00029-010-0017-z
  [arXiv:0906.0620 [math.QA]].
\bibitem{GK94}
  D.~Gepner and A.~Kapustin,
  ``On the classification of fusion rings,''
  Phys. Lett. B \textbf{349}, 71-75 (1995)
  doi:10.1016/0370-2693(95)00172-H
  [arXiv:hep-th/9410089 [hep-th]].
\bibitem{BNRW15}
  P.~Bruillard, S.H.~Ng, E.C.~Rowell, and Z.~Wang, ``ON CLASSIFICATION OF MODULAR CATEGORIES BY RANK,'' [arxiv:1507.05139 [math.QA]].
\bibitem{O00}
   A.~Ocneanu, ``The classification of subgroups of quantum $SU(N)$,'' Quantum Symmetries in Theoretical Physics and Mathematics R. Coquereaux et al (ed), American Mathematical Society, Providence, 2002, pp.133–159.
\bibitem{EP09}
   D.~E.~Evans and M.~Pugh, ``$SU(3)$-Goodman-de la Harpe-Jones subfactors and the realization of $SU(3)$ modular invariants,'' Rev. Math. Phys. 21 (2009), 877–928. https://doi.org/10.1142/S0129055X09003761
   [arXiv:0906.4252 [math.OA]].
\bibitem{CEM23}
  D.~Copeland and C.~Edie-Michell, ``CELL SYSTEMS FOR Rep(Uq(slN )) MODULE CATEGORIES,''
  [arXiv:2301.13172 [math.QA]].
\bibitem{G23}
  T.~Gannon,
  ``Exotic quantum subgroups and extensions of affine Lie algebra VOAs -- part I,''
  [arXiv:2301.07287 [math.QA]].
\bibitem{EM21}
  C.~Edie-Michell, ``TYPE II QUANTUM SUBGROUPS OF slN . I: SYMMETRIES OF LOCAL MODULES,'' Communications of the American Mathematical Society 3(2023), 112-165. https://doi.org/10.1090/cams/19
  [arXiv:2102.09065 [math.QA]].
\bibitem{FRS03}
  J.~Fuchs, I.~Runkel and C.~Schweigert,
  ``TFT construction of RCFT correlators. 2. Unoriented world sheets,''
  Nucl. Phys. B \textbf{678}, 511-637 (2004)
  doi:10.1016/j.nuclphysb.2003.11.026
  [arXiv:hep-th/0306164 [hep-th]].
\bibitem{KL02}
  Y.~Kawahigashi and R.~Longo,
  ``Classification of local conformal nets: Case c \ensuremath{<} 1,''
  Annals Math. \textbf{160}, 493-522 (2004)
  [arXiv:math-ph/0201015 [math-ph]].
\bibitem{FRS02}
  J.~Fuchs, I.~Runkel and C.~Schweigert,
  ``TFT construction of RCFT correlators 1. Partition functions,''
  Nucl. Phys. B \textbf{646}, 353-497 (2002). doi:10.1016/S0550-3213(02)00744-7
  [arXiv:hep-th/0204148 [hep-th]].
\bibitem{BBW21}
  D.~Barter, J.~C.~Bridgeman and R.~Wolf, ``Computing associators of endomorphism fusion categories,'' SciPost Phys. 13, 029 (2022) doi: 10.21468/SciPostPhys.13.2.029
  [arXiv:2110.03644 [math.QA]].
\bibitem{CZW18}
  S.~X.~Cui, M.~S.~Zini and Z.~Wang, ``On generalized symmetries and structure of modular categories,'' Sci China Math, 62(2019), 417–446. https://doi.org/10.1007/s11425-018-9455-5
  [arXiv:1809.00245 [math.QA].
\bibitem{D02}
  P.~Deligne, ``Catégories tensorielles,'' Moscow Math. Journal 2(2002) , 227-248. https://doi.org/10.17323/1609-4514-2002-2-2-227-248
\bibitem{KK23MFCrank9}
  K.~Kikuchi,
  ``Classification of connected \'etale algebras in multiplicity-free modular fusion categories up to rank nine,''
  [arXiv:2404.16125 [math.QA]].
\bibitem{TW19}
  R.~Thorngren and Y.~Wang,
  ``Fusion Category Symmetry I: Anomaly In-Flow and Gapped Phases,''
  [arXiv:1912.02817 [hep-th]].
\bibitem{HLS21}
  T.~C.~Huang, Y.~H.~Lin and S.~Seifnashri,
  ``Construction of two-dimensional topological field theories with non-invertible symmetries,''
  JHEP \textbf{12}, 028 (2021)
  doi:10.1007/JHEP12(2021)028
  [arXiv:2110.02958 [hep-th]].
\bibitem{O05}
  V.~Ostrik, ``Pre-modular categories of rank 3,'' Mosc. Math. J., 8(2008), 111–118. https://doi.org/10.17323/1609-4514-2008-8-1-111-118
  [arXiv:math/0503564 [math.QA]].
\bibitem{V88}
  C.~Vafa, ``Toward classification of conformal theories,'' Phys.Lett.B 206 (1988) 421–426. https://doi.org/10.1016/0370-2693(88)91603-6
\end{thebibliography}
\end{document}